\documentclass{gtart}

\usepackage{amsmath,amsfonts,amscd,flafter,rlepsf}

\newtheorem{theorem}{Theorem}[section]
\newtheorem{prop}[theorem]{Proposition}

\newtheorem{lemma}[theorem]{Lemma}

\theoremstyle{remark}

\newtheorem{remark}[theorem]{Remark}

\newcommand{\Z}{\mathbb{Z}}

\newcommand{\cm}{\cdot}

\newcommand{\ModSWfour}{\mathcal{M}}
\newcommand{\ModFlow}{\ModSWfour}

\newcommand{\SpinC}{{\mathrm{Spin}}^c}

\newcommand\Hom{\mathrm{Hom}}

\newcommand\abuts\Rightarrow
\newcommand\Sym{\mathrm{Sym}}

\newcommand{\F}{\mathbb F}

\newcommand{\ev}{\mathrm{ev}}
\newcommand{\odd}{\mathrm{odd}}

\newcommand\ModSphere{\ModFlow\left({\mathbb S}\longrightarrow 
\Sym^{g-1}(\Sigma_{1})\times \Sym^2(\Sigma_{2})\right)}
\newcommand\ModSpheres\ModSphere

\newcommand\HFp{\HFb}

\newcommand\HFinf{HF^\infty}

\newcommand\HFb{HF^+}

\newcommand\UnparModSp{\widehat \ModSp}
\newcommand\UnparModFlow\UnparModSp
\newcommand\Mod\ModSp

\newcommand\PD{\mathrm{PD}}

\newcommand{\spinc}{\mathfrak s}

\newcommand{\spinct}{\mathfrak t}

\newcommand\ModMaps{\mathcal M}
\newcommand\ModSp\ModMaps

\newcommand\Fp[1]{F^{+}_{#1}}

\newcommand\Field{\mathbb F}

\newcommand\Dual{\mathcal D}
\newcommand\Duality\Dual

\newcommand\Tor{\mathrm{Tor}}

\newcommand\Vertices{\mathrm{Vert}}
\newcommand\Char{\mathrm{Char}}
\newcommand\Combp{{\mathbb H}^+}
\newcommand\DCombp{{\mathbb K}^+}

\newcommand\NatTransp{T^+}
\newcommand\InjMod[1]{{\mathcal T}^+_{#1}}

\newcommand\MapOnep{A^+}
\newcommand\MapTwop{B^+}
\newcommand\MapThreep{C^+}

\newcommand\MapOneCombp{{\mathbb A}^+}
\newcommand\MapTwoCombp{{\mathbb B}^+}

\begin{document}

\title{On Heegard Floer homology of plumbed \\ three-manifolds with $b_1=1$}

\author{Raif Rustamov}
\address{The Program in Applied and Computational Mathematics\\Princeton University,New Jersey 08540, USA}
\email{rustamov@princeton.edu}

\begin{abstract}
In this paper we extend the results of \cite{Plumbing} to calculate Ozsv\'ath-Szab\'o Floer homology group $\HFp$ for a class of negative-semidefinite plumbings with $b_1 = 1$.
\end{abstract}
\maketitle

\section {Introduction}
In \cite{HolDisk} and \cite{HolDiskFour} topological invariants for closed oriented three manifolds and cobordisms between them are defined by using a construction from symplectic geometry. The resulting Floer homology package has many properties of a topological quantum field theory.

The construction of Heegard Floer homology is more combinatorial in flavor than the corresponding gauge theoretical constructions of Donaldson-Floer (see~\cite{DonaldsonFloer}) and Seiberg-Witten theories (see~\cite{Witten}, \cite{Morgan}, \cite{KMthom}).
However, the construction still depends on profoundly analytic objects - holomorphic disks. As a result, one is tempted to consider classes of manifolds which allow a completely combinatorial description of their Heegaard Floer homologies.

In \cite{Plumbing} a class of plumbing three-manifolds is studied. It is proved there that $\HFp$ of these manifolds can be expressed in terms of equivariant maps. We extend the results of this paper to get a similar description for negative-semidefinite plumbed manifolds with one bad vertex and $b_1=1$. To formulate the main theorems we first give required preliminaries.

Let G be a weighted graph and let $m(v)$ and $d(v)$ be respectively the weight and the degree of the vertex $v$.
We denote by $X(G)$ the four-manifold with boundary having $G$ as its plumbing diagram. Let $Y(G)$ be the oriented three-manifold which is the boundary of $X(G)$.

For $X=X(G)$, the group $H_2(X;\Z)$ is the lattice freely spanned by vertices of $G$.
 Denoting by $[v]$ the homology class in $H_2(X;\Z)$ corresponding to the vertex $v$ of $G$, the values of the intersection form of $X$ on the basis are given by $[v]\cm[v] = m(v)$; 
 $[v]\cm[w] = 1$ if $vw$ is an edge of $G$ and $[v]\cm[w]=0$ otherwise. The graph $G$ is called \emph{negative-(semi)definite} if the form is negative-(semi)definite. A vertex $v$ is said to be a \emph{bad vertex} of $G$ if $m(v)>-d(v)$.

Let $\InjMod{0}$ denote the $\Z[U]$-module which is the quotient of
$\Z[U,U^{-1}]$ by the submodule $U\cm \Z[U]$, graded so that the element $U^{-d}$ (for $d\geq 0$) is supported in degree $2d$. Define $\InjMod{d}=\InjMod{0}[d]$. 

Denoting by $\Char(G)$  the set of characteristic vectors for the
intersection form define 
$$\Combp(G)\subset \Hom(\Char(G),\InjMod{0})$$
to be the set of finitely supported functions satisfying the following relations for all characteristic vectors $K$ and vertices $v$:

\begin{equation}
\label{eq:AdjRel}
U^n\cm \phi(K+2\PD[v]) = \phi(K),
\end{equation}
if $2n=\langle K,v \rangle + v\cm v \geq 0$ and 
\begin{equation}
\label{eq:AdjRel2}
\phi(K+2\PD[v]) = U^{-n}\cm \phi(K)
\end{equation}
for $n<0$.

We can decompose $\Combp(G)$ according to $\SpinC$ structures over $Y$.
Note first that the first Chern class gives an identification of the set of
$\SpinC$ structures over $X=X(G)$ with the set of characteristic vectors
$\Char(G)$. Observe that the image of $H^2(X,\partial X;\Z)$ in $H^2(W;\Z)$
is spanned by the Poincar\'e duals of the spheres corresponding to the vertices.
Using the restriction to boundary, it is easy to see that the set of $\SpinC$ structures
over $Y$ is identified with the set of $2 H^2(X,\partial X;\Z)$-orbits 
in $\Char(G)$.

Fix a $\SpinC$ structure $\spinct$ over $Y$. Let $\Char_\spinct(G)$
denote the set of characteristic vectors for $X$ which are first
Chern classes of $\SpinC$ structures $\spinc$ whose restriction to the
boundary is $\spinct$. Similarly, we let $$\Combp(G,\spinct)\subset
\Combp(G)$$ be the subset of maps which are supported on the subset of
characteristic vectors $\Char_\spinct(G)\subset \Char(G)$.  We have a
direct sum splitting: $$\Combp(G)\cong \bigoplus_{\spinct\in\SpinC(Y)}
\Combp(G,\spinct).$$

For a negative definite graph $G$ a grading on $\Combp(G)$ is introduced as follows: we say that an element $\phi\in \Combp(G)$ is homogeneous of degree $d$ if 
for each characteristic vector $K$ with $\phi(K)\neq 0$, $\phi(K)\in \InjMod{0}$ is a homogeneous element with:
\begin{equation}
\label{eq:DefOfDegree}
\deg(\phi(K))-\frac{K^2+|G|}{4}=d.
\end{equation}

The following theorem is proved in \cite{Plumbing}:
\begin{theorem}
\label{intro:SomePlumbings}
Let $G$ be a negative-definite weighted forest with at most one bad vertex. Then, for each $\SpinC$ structure
$\spinct$ over $-Y(G)$, there is an isomorphism of graded
$\Z[U]$ modules,
$$\HFp(-Y(G),\spinct)\cong \Combp(G,\spinct).$$
\end{theorem}

We formulate a similar statement for negative-semidefinite trees with one bad point. In this case $b_1=1$ and so we have both torsion and non-torsion $\SpinC$ structures on $Y$. Note that for $K \in \Char_\spinct(G)$ with $\spinct$ torsion, $K^2$ is well defined. As a result, we can introduce a grading on $\Combp(G,\spinct)$ as before with \eqref{eq:DefOfDegree} replaced by 
\begin{equation}
\label{eq:DefOfDegreeSemidef}
\deg(\phi(K))-\frac{K^2+|G|-3}{4}=d.
\end{equation}

For a torsion $\SpinC$ structure Ozsv\'ath-Szab\'o Floer homology groups of manifolds with $b_1=1$ have absolute grading by half integers. This grading lifts absolute $\Z_2$ grading of the same groups: $-1/2$ modulo $2$ of the first correspond to the even degree of the latter and $1/2$ modulo $2$ corresspond to the odd degree.

The following theorem is an analogue of theorem~\ref{intro:SomePlumbings} for the negative-semidefinite case.

\begin{theorem}
\label{intro:Main}
Let $G$ be a negative-semidefinite weighted forest with at most one bad vertex. Then, for each $\SpinC$ structure
$\spinct$ over $-Y(G)$, there is an isomorphism of graded
$\Z[U]$ modules,
$$\HFp_{\odd}(-Y(G),\spinct)\cong \Combp(G,\spinct).$$

For non-torsion $\SpinC$ structure $\spinct$ we have 
$$\HFp_{\ev}(-Y(G),\spinct)\cong 0,$$
while if $\spinct$ is torsion then
$$\HFp_{\ev}(-Y(G),\spinct) \cong \InjMod{d}$$
for some $d=d(\spinct)$.
\end{theorem}

Let $\spinct$ be a torsion $\SpinC$ structure on a three-manifold $Y$ with $H_1(Y,\Z) \cong \Z$. Correction terms $d_{\pm 1/2}(Y,\spinct)$ are defined as the minimal grading of any non-torsion element in the image of $\HFinf(Y,\spinct)$ in $\HFp(Y,\spinct)$ with grading $\pm 1/2$ modulo 2. For the type of manifolds we are considering we see that coditions of being non-torsion and lying in the image of $\HFinf$ are superficial. An important property of these correction terms is given by 
$$d_{\pm1/2}(Y,\spinct)=-d_{\mp1/2}(-Y,\spinct),$$
see \cite{AbsGraded}, Section 4.2.

It follows from theorem~\ref{intro:Main} that $$d_{-1/2}(-Y(G),\spinct)=d(\spinct)$$ and 
$$d_{1/2}(-Y(G),\spinct)= \min_{\{K\in\Char_\spinct(G)\}} -\frac{K^2+|G|-3}{4}.$$

The above result gives a practical calculation of $d_{1/2}(-Y(G),\spinct)$: for
a given $\spinct\in\SpinC(Y)$, it is easy to see that the minimum in question is always achieved among the finitely many
characteristic vectors $K\in\Char_\spinct(G)$ with $$|K\cm v| \leq
|m(v)|.$$ (A smaller set containing these minimal vectors is described
in Section 3.1.)
 
For our calculations of $\HFp$ to be completed we need to specify $d(\spinct)$ for each torsion $\SpinC$ structure $\spinct$. Sometimes one can use truncated Euler characteristic of $\HFp$, calculated from Turaev torsion, see \cite{HolDiskTwo}, Section 10.6. Another way is to consider a negative-semidefinite forest $H$ corresponding to the manifold $-X(G)$. If $H$ has only one bad vertex, then one can use properties of correction terms to finish the calculation. 

The above results are proved in Section 2. In Section 3 we discuss the calculation of $\Combp(G)$. We end this paper with several sample calculations.

\medskip
\noindent{\bf{Acknowledgments}}\qua No words can express my gratitude to my advisor  Zolt\'an Szab\'o who kindly and patiently taught almost everything I know on Heegaard Floer homology and beyond. I would like to thank Paul Seymour for his invaluable support and encouragement.
\section{Proof of Theorem~\ref{intro:Main}}
\subsection{Theorem~\ref{intro:Main} over $\Z/2\Z$}
Let us postpone consideration of $\Z$ coefficients and concentrate on the case of $\F=\Z/2\Z$ coefficients. Thus, unless stated otherwise, all Floer homology groups in this subsection are meant to be taken with $\F$ coefficients(which we supress from the notation). In particular,  $\InjMod{0}$ now denotes the quotient of $\F[U,U^{-1}]$ by the submodule $U\cm \F[U]$.

Define a map $$\NatTransp\colon \HFp(-Y(G))\longrightarrow
\Combp(G)$$ as follows. The plumbing diagram gives a 
 cobordism from $-Y(G)$ to $S^3$. Now let $$\NatTransp(\xi)\colon
\Char(G)\longrightarrow \InjMod{0}$$ be the map given by
$$\NatTransp(\xi)(K)=\Fp{W(G),\spinc}(\xi)\in\HFp(S^3)=\InjMod{0},$$ where
$\spinc\in\SpinC(W(G))$ is the $\SpinC$ structure whose first Chern
class is $K$, and $\Fp{W(G),\spinc}$ denotes the four-dimensional
cobodism invariant defined in~\cite{HolDiskFour}.

Now we state a more precise statement of the first part of therem~\ref{intro:Main}.

\begin{theorem}
\label{thm:Odd}
Let $G$ be a negative-semidefinite tree with at most one bad vertex. Then,
$T^+$ induces a grading-preserving isomorphism:
$$\HFp_{\odd}(-Y(G),\spinct) \cong \Combp(G,\spinct)$$
\end{theorem}

The following proposition is proved in the same way as in \cite{Plumbing}.

\begin{prop}
\label{prop:SetUp}
The map $\NatTransp$ induces an $\Field[U]$-equivariant, degree-preserving
map from $\HFp_{\odd}(-Y(G),\spinct)$ to $\Hom(\Char_{\spinct}(G),\InjMod{0})$ whose image lies
in $$\Combp(G,\spinct)\subset \Hom(\Char_\spinct(G),\InjMod{0}).$$
\end{prop}

If $G$ is a weighted graph with a distinguished vertex
$v\in\Vertices(G)$, let $G'(v)$ be a new graph formed by
introducing one new vertex $e$ labelled with weight $-1$, and
connected to only one other vertex, $v$. Let $G_{+1}(v)$
denote the weighted graph whose underlying graph agrees with $G$, but
whose weight at $v$ is increased by one (and the weight stays the same
for all other vertices).
The two three-manifolds $Y(G'(v))$ and $Y(G_{+1}(v))$ are
diffeomorphic.

The following long exact sequence follows from Theorem~\ref{HolDiskTwo:thm:GeneralSurgery} of~\cite{HolDiskTwo}
\begin{gather*}
\cdots \longrightarrow\HFp(-Y(G'(v)))
\stackrel\MapOnep{\longrightarrow}
\HFp(-Y(G))\hspace{1.5in}\\
\hspace{1in}\stackrel{\MapTwop}{\longrightarrow}
\HFp(-Y(G-v))
\stackrel{\MapThreep}{\longrightarrow}
\HFp(-Y(G'(v))) \longrightarrow \cdots
\end{gather*}
Here the maps $\MapOnep$, $\MapTwop$, and $\MapThreep$
are induced by two-handle additions.
Note that $\MapOnep$ changes $\Z_2$ grading, while the other two maps preserve it. 

Without loss of generality we can suppose that our negative-semidefinite forest with one bad vertex is $G_{+1}(v)$ for negative-definite forest $G$ with (at most) one bad vertex. Corollary 1.4 of \cite{Plumbing} states that for negative-definite trees with at most one bad vertex $\HFp$ is supported in even degrees, and so , for example, $\HFp_{\odd}(-Y(G))$ is 0. From previous exact sequence we get
\begin{equation}
\label{eq:exactseq}
\begin{CD}
0 \longrightarrow\HFp_{\odd}(-Y(G'(v)))
\stackrel\MapOnep{\longrightarrow}
\HFp_{\ev}(-Y(G))\hspace{1.5in}\\
\hspace{1in}\stackrel{\MapTwop}{\longrightarrow}
\HFp_{\ev}(-Y(G-v))
\stackrel{\MapThreep}{\longrightarrow}
\HFp_{\ev}(-Y(G'(v))) \longrightarrow 0.
\end{CD}
\end{equation}

Corressponding to maps $\MapOnep$, $\MapTwop$  maps $\MapOneCombp$ and $\MapTwoCombp$ defined in \cite{Plumbing} so that the following diagram commutes:
\begin{equation}
\label{eq:CommDiagram}
\begin{CD}
\HFp_{\odd}(-Y(G'(v))) @>\MapOnep>>\HFp_{\ev}(-Y(G)) @>{\MapTwop}>> \HFp_{\ev}(-Y(G-v)),\\
@V{\NatTransp_{G'(v)}}VV @V{\NatTransp_{G}}VV
@V{\NatTransp_{G-v}}VV \\
\Combp(G'(v))) @>\MapOneCombp>>\Combp(G)
@>\MapTwoCombp>> \Combp (G-v).
\end{CD}
\end{equation}
Moreover $\MapOneCombp$ is injective, and $$\MapTwoCombp\circ\MapOneCombp=0.$$
The proof of these facts proceeds in the same way as the proof of Lemma 2.10 of \cite{Plumbing}. 

Gathering all of these together we have the diagram:
$$
\begin{CD}
0 @>>>\HFp_{\odd}(-Y(G'(v))) @>\MapOnep>>\HFp_{\ev}(-Y(G)) @>{\MapTwop}>> \HFp_{\ev}(-Y(G-v))\\
&& @V{\NatTransp_{G'(v)}}VV @V{\NatTransp_{G}}V{\cong}V
@V{\NatTransp_{G-v}}V{\cong}V \\
0@>>>  \Combp(G'(v))) @>\MapOneCombp>>\Combp(G)
@>\MapTwoCombp>> \Combp (G-v),
\end{CD}
$$

where the maps are indicated as isomorphisms when it follows from \cite{Plumbing}. From the diagram it follows that $T^+_{G'(v)}$ is an isomorphism. This concludes the proof of theorem~\ref{thm:Odd}.

To complete the proof of the main theorem we need to calculate $\HFp_{\ev}$. This will be done separately for the torsion and non-torsion case in following two lemmas.

\begin{lemma}
If $\spinct$ is a non-torsion $\SpinC$ structure then 
$$\HFp_{\ev}(-Y(G'(v)),\spinct)\cong 0.$$
\end{lemma}
\begin{proof}
Let $v$ be the bad vertex. We have a surjection $\MapThreep$ from  $\HFp_{\ev}(-Y(G-v))$ to $\HFp_{\ev}(-Y(G'(v))$. Since $G-v$ does not have any bad vertices its homology is given by the direct sum of $\InjMod{d}$'s, with one summand with its own $d$ for each $\SpinC$ structure on $-Y(G-v)$. This clearly means that for any integer $n$ and an element $\eta \in \HFp(-Y(G-v))$ there exists another element $\zeta$ so that $\eta=U^n \cdot \zeta$. Since $\spinct$ is non-torsion, then there exists an integer $N$ such that action of $U^N$ on  $\HFp_{\ev}(-Y(G'(v)),\spinct)$ is zero. Taking into account that $\MapThreep$ is $U$-equivariant and surjective proves the lemma.
\end{proof}
\begin{lemma}
If $\spinct$ is torsion then
$$\HFp_{\ev}(-Y(G'(v)),\spinct) \cong \InjMod{d}$$
for some integer d.
\end{lemma}
\begin{proof}
Given $b_1=1$ we know the group $\HFinf$ and so here it is enough to show that for any integer $n$ and an element $\eta \in \HFp(-Y(G'(v))$ there exists another element $\zeta$ so that $\eta=U^n \cdot \zeta.$ This again follows from the equivariance and surjectivity of $\MapThreep$.
\end{proof}
\subsection{Theorem~\ref{intro:Main} over $\Z$}
Since the sequence~\ref{eq:exactseq} holds with $\Z$ coefficients then the previous two lemmas are also true over $\Z$, with the proofs being the same. It remains to figure out $\HFp_{\odd}$. 
\begin{lemma}
$\HFp_{\odd}(-Y(G'(v)))$ over $\Z$ is torsion free.
\end{lemma}
\begin{proof}
The lemma easily follows from the sequence~\ref{eq:exactseq} over $\Z$ and the fact that $\HFp_{\ev}(-Y(G))$ is torsion free, see \cite{Plumbing}.
\end{proof}

\section{Calculations}
\subsection{Computing $\Combp(G)$}
 For calculational purposes it is helpful to adopt a dual point of view. Let $\DCombp(G)$ be the set of equivalence
 classes of elements of $\Z^{\geq 0} \times \Combp(G)$ (and we write $U^m\otimes K$ for the pair $(m, K)$) 
under the following equivalence relation.
For any vertex $v$ let
$$2n=\langle K,v \rangle + v\cm v.$$
If $n\geq 0$, then
\begin{equation}
\label{equation:rel1}
U^{n+m}\otimes (K+2\PD[v]) \sim U^m\otimes K,
\end{equation}
while if $n\leq 0$, then
\begin{equation}
\label{equation:rel2}
U^m\otimes (K+2\PD[v]) \sim U^{m-n}\otimes K.
\end{equation}

Starting with a map $$\phi\colon \Char(G)\longrightarrow \InjMod{0},$$
consider an induced map $${\widetilde \phi}\colon \Z^{\geq 0}\times
\Char(G)\longrightarrow \InjMod{0}$$
defined by $${\widetilde \phi}(U^n\otimes K)=U^n\cm \phi(K).$$
Clearly, the set of finitely-supported functions $\phi\colon \Char(G)\longrightarrow
\InjMod{0}$ whose induced map ${\widetilde \phi}$ descends to
$\DCombp(G)$ is precisely $\Combp(G)$.

A \emph{basic element} or \emph{basic vector} of $\DCombp(G)$ is one whose equivalence class does not contain  any element of form $U^m\otimes K$ with $m>0$.
Given two non-equivalent basic elements $K_1 = U^0 \otimes K_1$ and $K_2 = U^0 \otimes K_2$ in the same $\SpinC$ structure, one can find positive
 integers $n$ and $m$ such that 
$$ U^{n} \otimes K_1 \sim U^{m} \otimes K_2.$$ 
If, moreover, the numbers $n$ and $m$ are minimal then this relation will be called the \emph{minimal relationship} between $K_1$ and $K_2$.
\begin{remark}
\label{rem:matrix}
Let $M$ be the incidence matrix of our graph, with the diagonal elements equal to the weights of the corresponding vertices. Note that $K$ is torsion if and only if it is in the range of $M$. As a result, it is possible to set (without ambiguity) $K^2 = K a $ for any $a$ satisfying  $Ma= K^T.$
Non-torsion vectors lie in the complement of the range of $M$. It is an exercise in  linear algebra to show that for any vector $K$ corressponding to a non-torsion $\SpinC$ structure there exist non-negative integers $n$ and $m$, $m>n$ such that $$ U^{n} \otimes K \sim U^{m} \otimes K.$$
Such  relations with minimal $m$ and $n$ for non-torsion basic vectors $K$ are also included among the minimal relationships; we define the \emph{height} of $K$ to be minimal such $n$. 
\end{remark}

On can see that $\DCombp(G)$ is specified as soon as one finds its basic elements and all of the minimal relationships.
We describe now a modification of the algorithm for calculating the basic elements given in \cite{Plumbing} to the case of negative-semidefinite trees.

First of all we find all of the torsion basic vectors.
Let $K$ corresspond to a torsion $\SpinC$ structure and satisfy 
\begin{equation}
\label{eq:PartBox}
m(v) \leq \langle K, v\rangle \leq -m(v).
\end{equation}

Construct a sequence of vectors  
$K=K_0,K_1,\ldots,K_n$, where $K_{i+1}$ is obtained from $K_i$ by choosing any vertex $v_{i+1}$ with
$$\langle K_i,v_{i+1}\rangle = -m(v_{i+1}),$$
and then letting
$$K_{i+1}=K_i+2\PD[v_{i+1}].$$
Note that any two vectors in this sequence are equivalent.

Either the sequence is infinite or it terminates in one of two ways:
\begin{itemize}
\item the final vector $L=K_n$ satisfies the inequality,
\begin{equation}
\label{eq:OtherPartBox}
m(v) \leq \langle L, v\rangle \leq -m(v)-2
\end{equation}
at each vertex $v$ 
\item there is some vertex $v$ for which
\begin{equation}
\label{eq:OtherTermination}
\langle K_{n},v \rangle > -m(v).
\end{equation}
\end{itemize}
\begin{remark}
If the sequence is infinite then it can be made periodic. One just have to choose an ordering on the vertices of $G$. The ordering is used to eliminate the indeterminacy when there are several choices for $v_i$ while constructing the sequence.  
\end{remark}
Torsion vector $K$ is called \emph{good} if it satisfies inequality~\eqref{eq:PartBox} and either
\begin{itemize}
\item $\langle K, v\rangle \neq m(v)$ for all $v$ and previous sequence terminates in a characteristic vector $L$ satisfying
inequality~\eqref{eq:OtherPartBox}, or
\item there is $v$ with $\langle K, v\rangle = m(v)$ and the sequence is periodic.
\end{itemize}

\begin{prop} The equivalence classes in $\bigoplus_{c_1(\spinct) \in \Tor}
\Combp(G,\spinct)$ which
have no representative of the form $U^m\otimes K'$ with $m>0$ are in
one-to-one correspondence with good torsion vectors.
\end{prop}
\begin{proof}
The proof is similar to the proof of Proposition 3.2 of \cite{Plumbing}. The only new ingredient is that the sequence can become periodic because of the tree being semidefinite. In this case, since one can also subtract $\PD$'s instead of adding them, it follows that the sequence is purely periodic. As a consequence, the initial vector is basic.
\end{proof}
To find the non-torsion basic elements we consider all $K$ coressponding to non-torsion $\SpinC$ structures with 
\begin{equation}
\label{eq:PartBoxNon}
m(v)+2 \leq \langle K, v\rangle \leq -m(v).
\end{equation}
Then construct a sequence of vectors in the same way as for torsion case. This time the sequence cannot be periodic and it will terminate in one of the two ways described above.

\begin{prop} The equivalence classes in $\bigoplus_{c_1(\spinct) \in H^2-\Tor}
\Combp(G,\spinct)$ which
have no representative of the form $U^m\otimes K'$ with $m>0$ are in
one-to-one correspondence with non-torsion vectors $K$ satisfying
inequality~\eqref{eq:PartBoxNon} for which the previous sequence terminates in a characteristic vector $L$ satisfying
inequality~\eqref{eq:OtherPartBox} and for which there is no positive integer $l$ with $$K \sim U^l\otimes K.$$
\end{prop}

\begin{proof}
The proof is similar to the proof of Proposition 3.2 of \cite{Plumbing}. However, since $K^2$ cannot be defined in this case, the last part of the cited proof cannot be carried out. As a result, to get a one-one correspondence one has to exclude initial vectors $K$ which satisfy $K \sim U^l\otimes K$ for some positive integer $l$.  
\end{proof}

\subsection{Examples}
\begin{figure}[ht!]
\cl{\epsfxsize 1.0in\epsfbox{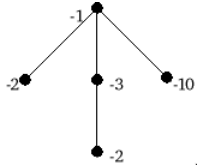}}
\caption{\label{fig:Example1}
{\bf{Plumbing description of 0 surgery on $(2,5)$ torus knot}}}
\end{figure}

{\bf{0-surgery on $(2,2n+1)$ torus knot}}\qua We give here another calculation of $\HFp(-Y_n)$ for a class of negative-semidefinite plumbed manifolds $Y_n$, where $n$ is a positive integer and $Y_n$ is obtained as 0-surgery on $(2,2n+1)$ torus knot, compare \cite{AbsGraded}. The plumbing diagram of $Y_2$ is depicted in figure~\ref{fig:Example1}. The diagram for $Y_n$ is similar, except it has $n-1$ of $-2$'s in the middle strand, and instead of $-10$ we have $4n-2$. 

We will write the
elements of $\Combp(G)$ as row vectors with the
first four coordinates corresponding to the vertices with weights $-1$, $-2$, $-3$ and $-4n-2$ respectively, and all remaining entries corresponding to $-2$'s on the middle strand ordered by the distance from the root starting with the closest one.

Let $E=E_n$ be a generator of $H^2(Y_n; \Z) = \Z$ and let $\spinct_i$ be the $\SpinC$ structure with first Chern class equal to $2iE$.  

Moreover, suppose that $$\Delta(T) =a_0+\sum_{i=1}^{d} a_i(T^i+T^{-i})$$ is the symmetrized Alexander polynomial of $(2,2n+1)$ torus knot. Define 
$$t_i=\sum_{j=1}^{d}j a_{|i|+j}.$$
One can easily calculate $t_i$ in the obvious way.

Let us start eith  $\HFp$ for $Y_2$. One calculates the basic vectors to be 
$$K_{-1}=(1,0,-1,-8,0),$$
$$K_0 =(1,0,-1,-6,0),$$
$$K_{1} = (1, 0, -1, -4,0).$$
$K_0$ is the only torsion vector, while  $K_{-1}$ and $K_1$ correspond to $\pm 2E$. One can show that $U \otimes K_{\pm 1} \sim U^2 \otimes K_{\pm 1}$, i.e. the height of $K_{\pm 1}$ is equal to 1. 
It follows at once that  $$\HFp(-Y_2, \spinct_{\pm1}) = \Z[U]/U^1,$$ supported in the odd degree, and
$$\HFp(-Y_2, \spinct) =0 $$
for all other non-torsion $\spinct$.
We have 
$$\HFp_{\odd}(-Y_2,\spinct_0) = \InjMod{1/2}.$$
From the second part of the main theorem it follows that $$\HFp_{\ev}(-Y_2,\spinct_0) = \InjMod{d}$$
for some $d$. 
Since $Y_2$ is a 0-surgery on (2,5) torus knot, we can easily calculate the truncated Euler characteristic of $\HFp(-Y_2,\spinct_0)$ from the Alexander polynomial, see~\cite{AbsGraded}. In fact, $\chi^{\mathrm{trunc}}=-t_0=-1$ and so $d=3/2$. 

Now let us turn to the general situation. Firstly, let us calculate all the basic vectors.

\begin{lemma}
\label{lemma:bv}
For $Y_n$ there are $2n-1$ non-equivalent basic vectors $K_{-n+1},K_{-n+2},$ $..., K_0, K_1, K_{n-1}$ where
$$K_i = (1, 0, -1, -2n-2+2i, 0, 0, ..., 0).$$
There is only one non-torsion basic vector - $K_0$.  
\end{lemma}
\begin{proof}
The proof is similar to the proof of lemma 3.1 of \cite{StupidPaper}. For the second part one has to find out the vectors lying in the range of $M$, see remark ~\ref{rem:matrix}.
\end{proof}
Now we should find all of the minimal relations.\begin{lemma}
\label{lemma:rel}
The following relationships are satisfied for non-torsion basic vectors:
$$U^{t_i} \otimes K_{-n+i} \sim U^{n+1-i} \otimes K_{-n+i},$$
$$U^{t_i} \otimes K_{n-i} \sim U^{n+1-i} \otimes K_{n-i},$$
where $i=1,2,...,n-1$. These relationships are minimal, i.e. height of $K_{\pm(n-i)}$ is equal to $t_i$.
\end{lemma}
\begin{proof}
The proof is similar to the proof of lemma 3.2 of \cite{StupidPaper}. The minimality follows from 
$$\chi(\HFp(-Y_n,\spinct_i))=-t_i$$ 
and the fact that all of $\HFp$ for non-torsion $\SpinC$ structures is supported in the odd degree.
\end{proof}
\begin{lemma}
\label{lemma:length}
The renormalized length $(K_0^2+|G|-3)/4$ of $K_0$ is equal $-1/2$.
\end{lemma}
\begin{proof}
One proves by induction that $K_0^2=-n-2$.
\end{proof}

We get the following result about the Floer homology groups of $-Y_n$:
\begin{prop}$$\HFp(-Y_n, \spinct_i)=\Z[U]/U^{t_i}$$
$$\HFp(-Y_n, \spinct_0) = \InjMod{1/2}\oplus \InjMod{2t_0-1/2}$$
\end{prop}
\begin{proof}
We only need to calculate $d$ for the even part of $\HFp$, and this is done using Euler characteristic argument as we did in the case of $-Y_2$.
\end{proof}

\begin{figure}[ht!]
\cl{\epsfxsize 2.0in\epsfbox{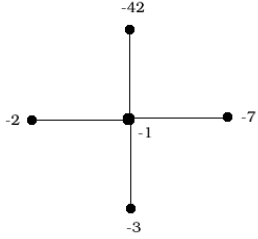}}
\caption{\label{fig:Example2}
{\bf{Plumbing diagram G}}}
\end{figure}

{\bf{Another example}}\qua Let us consider the plumbing diagram depicted in the figure  \ref{fig:Example2}. We will use shorthand $Y=Y(G)$. Note that there is another negative-semidefinite forest $H$ with one bad vertex satisfying  $X(G)=-X(H)$. Picture of $H$ is also star shaped, with $-3$ in the center, and four strands coming out of the center, first strand contains one $-2$, second - two $-2$'s, third - six $-2$'s and the last one with forty one $-2$'s. We will use $H$ to calculate $d$. 

We will again write the
elements of $\Combp(G)$ as row vectors with the
coordinates corresponding to the vertices with weights $-1$, $-2$, $-3$, $7$ and $-42$ in the order given. 
Let $E=E_n$ be a generator of $H^2(Y_n; \Z) = \Z$ and let $\spinct_i$ be the $\SpinC$ structure with first Chern class equal to $2iE$.  

Basic vectors are $$S_i=(1, 0,-1, -3 ,-10+2i)$$ and $$T_j =(1, 0, -1, -5, 2+2j)$$ for $i=-15,-14,..., 21$ and $j= -21,-25,...,15$. Moreover, $S_i$ lies in the $\SpinC$ structure $\spinct_i$, and $T_j$ lies in $\spinct_j$. Note that since everything is symmetric under $\spinct_i$ goes to $\spinct_{-i}$, we need to consider only $i \geq 0$. 

Let $h_i=1,5,4,4,3,3,3,3,2,2,1,1,1,...$ for $i=0,1,...21$. It turns out that heights of $S_i$ and $T_i$ are both equal to $h_i$ for $i \geq1$.
The minimal relations are given by 
$$ U^{h_i} \otimes S_{i} \sim U^{h_i} \otimes T_{i}$$
for $0\leq i \leq 15$. One calculates $S_0^2=T_0^2 = -4$.
We conclude that
$$ \HFp(-Y, \spinct_i)=0, \mathrm{for \qua} i \geq 22,$$
$$ \HFp(-Y, \spinct_i)=\Z[U]/U^{h_i},\mathrm{for \qua} 16\leq i \leq 21,$$
$$ \HFp(-Y, \spinct_i)=\Z[U]/U^{h_i} \oplus \Z[U]/U^{h_i},\mathrm{for \qua} 1 \leq i \leq 15, $$
$$\HFp_{\odd}(-Y, \spinct_0) = \InjMod{1/2} \oplus \Z_{1/2}.$$

To complete the calculation we need $d(\spinct_0)=d_{-1/2}(-Y, \spinct_0)$. Since $d_{-1/2}(-Y, \spinct_0)=-d_{1/2}(Y,\spinct_0)=-d_{1/2}(-Y(H),\spinct_0)$ we have that
$$d(\spinct_0)= - \min_{\{K\in\Char_{\spinct_0}(H)\}} -\frac{K^2+|H|-3}{4}=23/2.$$
We conclude 
$$\HFp(-Y, \spinct_0) = \InjMod{1/2} \oplus \InjMod{23/2} \oplus \Z_{1/2}.$$

\end{document}